\documentclass[12pt,amsmath]{amsart}
\usepackage{amssymb}
\theoremstyle{plain}
\newtheorem{theorem}{Theorem}[section]

\newtheorem{lemma}[theorem]{Lemma}

\newtheorem{proposition}[theorem]{Proposition}
\newtheorem{fact}[theorem]{Fact}
\newtheorem{claim}[theorem]{Claim}
\theoremstyle{definition}
\newtheorem{definition}[theorem]{Definition}
\newtheorem{remark}[theorem]{Remark}
\newtheorem{notation}[theorem]{Notation}
\newtheorem{convention}[theorem]{Convention}
\theoremstyle{remark}

\newcommand{\wilog}{without loss of generality}

\begin{document}

\title{On cardinalities in quotients of inverse limits of groups}


\author{Saharon Shelah}
\email[Saharon Shelah]{ shelah@sunset.ma.huji.ac.il}  
\address{Institute of Mathematics\\
 The Hebrew University\\
 Jerusalem, Israel\\
\&
 Rutgers University\\
 Mathematics Department\\
 New Brunswick, NJ 08903  USA}

\author{Rami Grossberg}\thanks{
Partially supported by the United States-Israel Binational Science
Foundation and by the National Science Fundation, NSF-DMS97-04477.}
\email[Rami Grossberg]{ rami@cmu.edu}

\address{
Department of Mathematical Sciences\\
Carnegie Mellon University\\
Pittsburgh, PA 15213}


\date{\today}

\begin{abstract} Let $\lambda$ be $\aleph_0$ or a strong limit of cofinality $\aleph_0$.
Suppose that $\langle G_m,\pi_{m,n}\;:\;m\leq n<\omega\rangle$ and 
 $\langle H_m,\pi^t_{m,n}\;:\;m\leq n<\omega\rangle$ are projective
 systems of groups of cardinality less than $\lambda$ and suppose that 
for every $n<\omega$ there is a homorphism $\sigma:H_n\rightarrow G_n$ such
that all the diagrams commute.

If for every $\mu<\lambda$ there exists
$\langle f_i\in G_{\omega} \;:\;i<\mu\rangle$ such that
 $i\neq j\Longrightarrow f_if_j^{-1}\not\in\sigma_{\omega}(H_{\omega})$
then there exists $\langle f_i\in G_{\omega} \;:\;i<2^{\lambda}\rangle$
such that $i\neq j\Longrightarrow f_if_j^{-1}\not\in\sigma_{\omega}(H_{\omega})$.

\end{abstract}
\maketitle

\section{Introduction}

The main result of this paper was motivated by our interest in the structure 
of the group $Ext_p(G,\mathbf {Z})$ for $G$ abelian torsion free.
For basic results about the structure of $Ext(G,\mathbf {Z})$ the reader
is refered to sections 47 and 52 of Laszlo  Fuchs book \cite{Fu}, however all we need is
Definition \ref{extp} below.  Since Shelah's proof of the independence of
Whitehead's problem of $ZFC$ (see \cite{Sh:44}) much was done since that paper,
for a summary
see the introduction to \cite{GrSh} and Chapter XII of Eklof \& Mekler's book
is dedicated (\cite{EK}) to the structure of $Ext$.

In \cite{GrSh} we have dealt with the cardinality of $Ext_p(G,\mathbf {Z})$.
The main Theorem of \cite{GrSh} states that for a strong limit $\lambda$
of cofinality $\aleph_0$ for every torsion free $G$ of cardinality $\lambda$
either 
\[
|Ext_p(G,\mathbf {Z})|<\lambda\text{ or } |Ext_p(G,\mathbf {Z})|=2^{\lambda}.
\]

In section 2 of \cite{GrSh} we  indicated  that the proof of the main
theorem can be adapted to give a result concerning cardinalities
of inverse systems of abelian groups subject to certain 
conditions (See Theorem \ref{1.1} below).  We did not include 
a proof there.  Recently we were asked to supply a complete proof
to that theorem.  Charles Megibben in a widely circulated
preprint \cite{meg} (which to our knowledge did not appear
yet in print) even claimed that he proved a result that contradicts
 Theorem \ref{1.1}.

The aim of this paper is to present a complete proof of Theorem \ref{1.1} below.

Notice that we do not make any assumptions on the groups, 
in particular the groups need not be commutative and can be even locally finite.
See more on the subject in \cite{Sh:664}.

\bigskip

\begin{theorem} \label{maintheorem}\label{1.1}[The Main Theorem] Suppose 
$\lambda$ is $\aleph_0$ or it is strong limit 
cardinal of cofinality $\aleph_0$.

\begin{enumerate}

\item     Let 
$\langle G_m,\pi_{m,n}:m \le n < \omega \rangle$ be an 
inverse system of groups of cardinality less than $\lambda$
whose inverse limit is $G_\omega$ with $\pi_{n,\omega}$ 
such that $|G_n| < \lambda$.  ( $\pi_{m,n}$ is
a homomorphism from $G_m$ to $G_n,\alpha \le \beta \le \gamma \le \omega
\Rightarrow \pi_{\alpha,\beta} \circ \pi_{\beta,\gamma} = 
\pi_{\alpha,\gamma}$ and $\pi_{\alpha,\alpha}$ is the identity).

\item   Let $\mathbf I$ be a finite index set.  Suppose that for every $t \in \mathbf I$,
 $\langle H^t_m,\pi^t_{m,n}:m \le n < \omega \rangle$ is an inverse system
of groups of cardinality less than $\lambda$ 
and $H^t_\omega$ with $\pi^t_{n,\omega}$ be the corresponding
inverse limit.

\item   Let for every $t \in \mathbf I,\sigma^t_n:H^t_n \rightarrow
G_n$ be a homomorphism such that all diagrams commute (i.e.
$\pi_{m,n} \circ \sigma^t_n = \sigma^t_m \circ \pi^t_{m,n}$ for $m \le n <
\omega$), and let $\sigma^t_\omega$ be the induced homomorphism from $H^t_
\omega$ into $G_\omega$.

\end{enumerate}

Assume that for every $\mu < \lambda$ there is a sequence
$\langle f_i \in G_\omega:i < \mu \rangle$ such that for 
$i \ne j $ and $t \in \mathbf I \Rightarrow f_if^{-1}_j \notin 
\text{ Rang}(\sigma^t_\omega)$. 
Then there is $\langle f_i \in G_\omega:i < 2^\lambda \rangle$ such
that 
$i \ne j $ and $t \in \mathbf I \Rightarrow f_i f^{-1}_j \notin 
\text{ Rang}(\sigma^t_\omega)$.

\end{theorem}
\bigskip

\begin{notation}  Since $\lambda$ has cofinality $\aleph_0$ we can fix
 $\lambda_n < \lambda$ for $n < \omega$ such
that $\lambda =  \sum_{n < \omega} \lambda_n$,  for all $n < \omega,
\lambda_n$ is regular and $2^{\lambda_n} < \lambda_{n+1} < \lambda$ and
$|G_n| +  \sum_{t \in \mathbf I} |H^t_n| \le \lambda_n$.

Denote by  $e_{G_\alpha},e_{H^t_\alpha}$  the unit elements.
Without loss of generality the groups are pairwise disjoint.
\end{notation}

\begin{definition}  
\begin{enumerate}

\item For $\alpha \le \omega$ let
$H_\alpha =  \prod_{t \in \mathbf I} H^t_\alpha$ and $H_{< \alpha} =
\prod_{\beta < \alpha} H_\beta$, 
$H_{\le \alpha} = \prod_{\beta \le \alpha}H_\beta$. 

\item For $\bar g \in H_\alpha$ let lev$(\bar g) = \alpha$, for $g \in 
H^t_\alpha$ let lev$(g) = \alpha$ (\wilog \, this is well defined). 

\item  For $\alpha \le \beta \le \omega,g \in H^t_\beta$ let $g \restriction
H^t_\alpha = \pi^t_{\alpha,\beta}(g)$ and we say \emph{$g \restriction H^t_\alpha$
is below $g$} and \emph{$g$ is above $g \restriction H^t_\alpha$} or \emph{extend
$g \restriction H^t_\alpha$}. 

\item For $\alpha \le \beta \le \omega,f \in G_\beta$ let $f \restriction
G_\alpha = \pi_{\alpha,\beta}(f)$.

\end{enumerate}
\end{definition}

We will now introduce the rank function used in the proof of Theorem \ref{1.1}, it 
is a measure for the possibility to extend functions in Lemma \ref{1.8} we show that it is
an ultrametric valuation.

\begin{definition} \label{1.5}
\begin{enumerate}

\item  For $g \in H^t_n,f \in G_\omega$ we
say that \emph{$(g,f)$ is a nice $t$-pair} if $\sigma^t_n(g) = f \restriction G_n$.

\item  Define a ranking function rk$_t(g,f)$ for any nice $t$-pair.  
First by induction on the ordinal $\alpha$ (we can fix $f \in G_\omega$), 
we define when rk$_t(g,f) \ge \alpha$ simultaneously for all $n < \omega$ and every 
$g \in H^t_n$

\begin{enumerate}

\item   rk$_t(g,f) \ge 0$ iff $(g,f)$ is a nice $t$-pair

\item   rk$_t(g,f) \ge \delta$ for a limit ordinal $\delta$ iff
for every $\beta < \delta$ we have rk$_t(g,f) \ge \beta$

\item   rk$_t(g,f) \ge \beta + 1$ iff 
$(g,f)$ is a nice $t$-pair,
and letting $n = \text{ lev}(g)$ there exists $g' \in H^t_{n+1}$ extending
$g$ such that rk$_t(g',f) \ge \beta$

\item   rk$_t(g,f) \ge -1$.

\end{enumerate}

\item For $\alpha$ an ordinal or $-1$ (stipulating $-1 < \alpha < \infty$ for
any ordinal $\alpha$) we have rk$_t(g,f) = \alpha$ iff
rk$_t(g,f) \ge \alpha$ and it is false that rk$_t(g,f) \ge \alpha +1$. 

\item rk$_t(g,f) = \infty$ iff for every ordinal $\alpha$ we have
rk$_t(g,f) \ge \alpha$. 

\end{enumerate}
\end{definition}

The following two claims give the principal properties of rk$_t(g,f)$.

\begin{claim}\label{1.6}  Let $(g,f)$ be a nice $t$-pair. 
\begin{enumerate}

\item The following statements are equivalent:
\begin{enumerate}

\item   rk$_t(g,f) = \infty$

\item   there exists $g' \in H^t_\omega$ extending $g$ such that
$\sigma^t_\omega(g') = f$.

\end{enumerate}

\item  If rk$_t(g,f) < \infty$, then rk$_t(g,f) < \lambda^+$. 

\item  If $g'$ is a proper extension of $g$ and $(g',f)$ is also a nice $t$-pair
then

\begin{enumerate}

\item  rk$_t(g',f) \le \text{ rk}_t(g,f)$ and

\item   if $0 \le \text{ rk}_t(g,f) < \infty$ then 
the inequality is strict.
\end{enumerate}

\end{enumerate}
\end{claim}

\begin{proof}
\begin{enumerate}

\item  $(a) \Rightarrow (b)$: 
Let $n$ be 
such that $g \in H^t_n$.  It is enough to define $g_k \in H^t_k$ for
$k < \omega,k \ge n$ such that  

\begin{enumerate}
\item   $g_n = g$

\item   $g_k$ is below $g_{k+1}$ that is $\pi^t_{k,k+1}
(g_{k+1}) = g_k$ and

\item   rk$_t(g_{k+1},f) = \infty$:
\end{enumerate}

Let $g' := \underrightarrow{\lim} g_k$ it is as required.  The 
definition is by induction on $k \ge n$. 
For $k = n$ let $g_0 = g$. 
For $k \ge n$, suppose $g_k$ is defined.  By (iii) we have rk$_t(g_k,f) =
\infty$, hence there exists $g^* \in H^t_{k+1}$ extending $g_k$ such that
rk$_t(g^*,f) = \infty$, and let $g_{k+1} := g^*$.

 $(b) \Rightarrow (a)$:  
Since $g$ is below $g'$, it is enough to prove by induction on $\alpha$ that
for every $k \ge n$ when $g_k := g' \restriction H^t_k$ we have that
rk$_t(g,f) \ge \alpha$.

For $\alpha = 0$, since $\sigma^t_\omega(g') = f \restriction G_n$ clearly for
every $k$ we have $\sigma^t_k(g_k) = f \restriction G_k$ so $(g_k,f)$ is a
nice $t$-pair.

For limit $\alpha$, by the induction hypothesis for every $\beta < \alpha$
and every $k$ we have rk$_t(g_k,f) \ge \beta$, hence by Definition 
\ref{1.5}(2)(b), rk$_t(g_k,f) \ge \alpha$.

For $\alpha = \beta +1$, by the induction hypothesis for every $k \ge n$ 
we have rk$_t(g_k,f) \ge \beta$.  Let $k_0 \ge n$ be given.  Since
$g_{k_0}$ is below $g_{k_0+1}$ and rk$_t(g_{k_0+1},f) \ge \beta$,
Definition \ref{1.5}(2)(c) implies that rk$_t(g_{k_0},f) \ge \beta +1$; i.e.
for every $k \ge n$ we have rk$_t(g_k,f) \ge \alpha$.  So we are done. 

\item
 Let $g \in H^t_n$ and $f \in G_\omega$ be given.  It is enough to prove
that if rk$_t(g,f) \ge \lambda^+$ then rk$_t(g,f) = \infty$.  Using part (1)
it is enough to find $g' \in H^t_\omega$ such that $g$ is below $g'$ and
$f = \sigma^t_\omega(g')$.

We define by induction on $k < \omega,g_k \in H^t_{n+k}$ such that
$g_k$ is below $g_{k+1}$ and rk$_t(g_k,f) \ge \lambda^+$.  For $k=0$ let
$g_k = g$.  For $k+1$, for every $\alpha < \lambda^+$, as rk$_t(g_k,f) >
\alpha$ by \ref{1.5}(2)(c) there is $g_{k,\alpha} \in G_{n+k+1}$ extending
$g_k$ such that rk$_t(g_{k,\alpha},f) \ge \alpha$.  But the number of possible
$g_{k,\alpha}$ is $\le |H^t_{n+k+1}| \le 2^{\lambda_{n +k+1}} < \lambda^+$
hence there are a function $g$ and a set $S \subseteq \lambda^+$ of 
cardinality $\lambda^+$ such that $\alpha \in S \Rightarrow g_{k,\alpha}
=g$.  Now take $g_{k+1} = g$.

\item Immediate.  
\end{enumerate}\end{proof}

\begin{lemma} \label{1.7} \begin{enumerate}

 \item Let $(g,f)$ be a nice $t$-pair.  Then
we have rk$(g,f) \le \text{ rk}(g^{-1},f^{-1})$. 

\item For every nice $t$-pair $(g,f)$ we have rk$(g,f) = \text{ rk}(g^{-1},
f^{-1})$.
\end{enumerate}
\end{lemma}

\begin{proof}   \begin{enumerate}

 \item By induction on $\alpha$ prove that rk$(g,f) \ge \alpha
\Rightarrow \text{ rk}(g^{-1},f^{-1}) \ge \alpha$ (see more details in
Lemma \ref{1.8}). 

\item  Apply part (1) twice. 
\end{enumerate}\end{proof}

In the following lemma we show that the rank is indeed  ultrametric
(ordinal valued).
 
\begin{lemma} \label{1.8}
 Let $n < \omega$ be fixed, and let $(g_1,f_1),
(g_2,f_2)$ be nice $t$-pairs with $g_\ell \in H^t_n(\ell = 1,2)$. 

\begin{enumerate}

\item  If $(g_1,f_1)$ and $(g_2,f_2)$ are $t$-nice pairs, then 
$(g_1g_2,f_1f_2)$ is a nice pair and rk$_t(g_1g_2,f_1f_2) \ge 
\text{ Min}\{\text{rk}_t(g_\ell,f_\ell):\ell =1,2\}$.

\item  Let $n,(f_1,g_1)$ and $(f_2,g_2)$ be as above.  If rk$_t(g_1,f_1) \ne
\text{ rk}_t(g_2,f_2)$, then 
rk$_t(g_1g_2,f_1f_2) = \text{ Min}
\{\text{rk}_t(g_\ell,f_\ell):\ell =1,2\}$.
\end{enumerate}
\end{lemma}

\begin{proof} 
\begin{enumerate}

\item It is easy to show that the pair is $t$-nice.  
We show by induction
on $\alpha$ simultaneously for all $n < \omega$ and every $g_1,g_2 \in H^t_n$
that Min$\{\text{rk}(g_\ell,f_\ell):\ell =1,2\} \ge \alpha$ implies that
rk$(g_1g_2,f_1f_2) \ge \alpha$.

When $\alpha = 0$ or $\alpha$ is a limit ordinal this is easy.  Suppose
$\alpha = \beta + 1$ and that rk$(g_\ell,f_\ell) \ge \beta +1$; by the
definition of rank for $\ell =1,2$ there exists $g'_\ell \in H^t_{n+1}$ 
extending $g_\ell$ such that $(g'_\ell,f_\ell)$ is a nice pair and 
rk$_t(g'_\ell,f_\ell) \ge \beta$.  By the induction assumption 
rk$_t(g'_1g'_2,f_1f_2) \ge \beta$.
Hence $g'_1g'_2$ is as required in the definition of rk$_t(g_1g_2,f_1f_2) 
\ge \beta + 1$. 

\item Suppose \, \wilog \, that rk$(g_1,f_1) < \text{ rk}(g_2,f_2)$, let
$\alpha_1 = \text{ rk}(g_1,f_1)$ and let $\alpha_2 = \text{ rk}_t(g_2,f_2)$.
By part (1), rk$_t(g_1g_2,f_1f_2) \ge \alpha_1$, by Proposition \ref{1.7},
rk$_t(g^{-1}_2,f^{-1}_2) = \alpha_2 > \alpha_1$.  So we have
\[
\alpha_1 = \text{ rk}_t(g_1,f_1) = \text{ rk}_t(g_1g_2g^{-1}_2,f_1f_2
f^{-1}_2) 
 \cr\ge \text{ Min}\{\text{rk}_t(g_1g_2,f_1f_2),\text{rk}_t(g^{-1}_2,
f^{-1}_2)\} 
  \cr = \text{ rk}_t(g_1g_2,f_1f_2) \ge \alpha_1.
\]
\end{enumerate}

Hence the conclusion follows. 
\end{proof}

\begin{definition} \begin{enumerate}\label{1.9}

\item  Let $\mu<\lambda$ and let 
$\bar \alpha = \langle \alpha_t:t \in \mathbf I \rangle$ where $\alpha_t$ is an
ordinal less or equal to $ \lambda^+$. We say that 
$\bar f = \langle f_i:i < \mu \rangle$ \emph{$\mu$-exemplifies $\bar \alpha
\in \Gamma_n$} (or \emph{$\bar f$ is a $\mu$-witness for $\bar \alpha \in \Gamma_n$})
iff
\begin{enumerate}

\item  $f_i \in G_\omega$ and $f_i \restriction G_n = e_{G_n}$

\item   for $i \ne j$ and $t \in \mathbf I$ we have rk$_t(e_{H^t_n},
f_if^{-1}_j) < \alpha_t$ (possibly is $-1$).

\end{enumerate}

\item Let
\[
\Gamma_n = \biggl\{ \bar \alpha\;:\;\bar \alpha = \langle \alpha_t:t \in
\mathbf I \rangle,\alpha_t \text{ an ordinal } \le \lambda^+, \cr
  \text{and for every } \mu < \lambda \text{ there is a sequence }
\langle f_i:i < \mu \rangle \cr
  \text{which } \mu \text{-exemplifies } \bar \alpha \in \Gamma_n
 \biggr\}.
\]

\item $\Delta_n = \{\bar \alpha \in \Gamma_n:\text{for no } \bar \beta
\text{ we have } \bar \beta \in \Gamma_n,\bar \beta \le \bar \alpha$ (i.e.
$ \bigwedge_{t \in \mathbf J_n} \beta_t \leq \alpha_t)$ and $\bar \beta \ne
\bar \alpha\}$.
\end{enumerate}
\end{definition}

\begin{claim}\label{1.10}  \begin{enumerate}

\item  $\Gamma_n$ is not empty.

\item $\Delta_n$ is not empty in fact $(\forall \bar \alpha \in \Gamma_n)
(\exists \bar \beta \in \Delta_n)(\bar \beta \le \bar \alpha)$.
\end{enumerate}
\end{claim}

\begin{proof}  \begin{enumerate}

\item  Let $\alpha^*_t = \sup\{\text{rk}_t(g,f)+1:g \in H^t_n,
f \in G^\omega \text{ and rk}_t(g,f) < \infty\}$, 
by \ref{1.6}(2), this is a supremum on a set of ordinals
$< \lambda^+$ (as $-1 +1=0$) hence is an ordinal $\le \lambda^+$.  So
$\langle \alpha^*_t:t \in \mathbf I \rangle$ is as required. 

\item  If not, then choose by induction on $\ell < \omega$ a sequence
$\bar \beta^\ell \in \Gamma_n$ such that $\bar \beta^0 = \bar \alpha,
\bar \beta^{\ell +1} \le \bar \beta^\ell,\bar \beta^{\ell +1} \ne \beta^\ell$.
So for each $t \in \mathbf I$, the sequence $\langle \beta^\ell_t:\ell <
\omega \rangle$ is a non-increasing sequence of ordinals hence is eventually
constant, say for some $\ell_t < \omega$ we have $\ell \in [\ell_t,\omega)
\Rightarrow \beta^\ell_t = \beta^{\ell_t}_t$, so as $\mathbf I$ is finite,
$\ell(*) = \text{ max}\{\ell_t:t \in \mathbf I\} < \omega$, so $\bar \beta
^{\ell(*)} = \bar \beta^{\ell(*)+1}$, a contradiction.
\end{enumerate}
\end{proof}

\begin{claim} \label{1.11}
\begin{enumerate}

\item  If $\mu \le \mu'$ and $\langle f_i:i < \mu'
\rangle,\mu'$-exemplify $\bar \alpha \in \Gamma_n$ and $h:\mu \rightarrow \mu'$
is one to one, then $\langle f_{h(i)}:i < \mu \rangle, \mu$-exemplifies
$\bar \alpha \in \Gamma_n$. 

\item   If $\langle f_i:i < \mu \rangle,\mu$-exemplify $\bar \alpha \in \Gamma_n$
and $f_i \restriction G_{n+1} = f$ for $i < \mu$, then $\langle f_i
f^{-1}_0:i < \mu \rangle,\mu$-exemplify $\bar \alpha \in \Gamma_{n+1}$. 

\item If $\bar \alpha \in \Gamma_n$, then $\bar \alpha \in \Gamma_{n+1}$.

\item If $\bar \alpha \in \Delta_n$, then some $\bar \beta \le \bar \alpha$
belongs to $\Delta_{n+1}$. 

\item For some $n < \omega$ there is $\bar \alpha \in \bigcap_{m \ge n} \Delta_n$.

\item In clause (b) of Definition \ref{1.9}(1) 
it suffices to deal with $i < j$.
\end{enumerate}
\end{claim}

\begin{proof}\begin{enumerate}
  
\item  Trivial. 

\item Clearly.

Clause (a):

$(f_i \circ f^{-1}_0) \restriction G_{n+1} = \sigma^\omega_{n+1}
(f_if^{-1}_0) = (\sigma^\omega_{n+1}(f_i))(\sigma^\omega_{n+1}(f_0))^{-1} = 
ff^{-1} = e_{G_{n+1}}$.

Clause (b):

 For $i \ne j$ and $t \in \mathbf I$, note that
\[
(f_if^{-1}_0)(f_jf^{-1}_0) = f_if^{-1}_0 f_0 f^{-1}_j =
f_i f^{-1}_j
\]
so we can use the assumption. 

\item   So let $\mu < \lambda$ and we should find a $\mu$-witness for
$\bar \alpha \in \Gamma_{n+1}$.  We can choose $\mu'$ such that $\mu \times
|G_{n+1}| < \mu' < \lambda$.  As $\bar \alpha \in \Gamma_n$, clearly there is
a $\mu'$-witness $\langle f_i:i < \mu' \rangle$ for it.  Now the number of
possible $f_i \restriction G_{n+1}$ is $\le |G_{n+1}|$ (really) even
$\le |\text{Rang}(\pi_{n+1,\omega}) \cap \text{ Ker}(\pi_{n,n+1})|)$ hence
for some $f \in G_{n+1}$ and $Y \subseteq \mu'$ we have: $|Y| \ge \mu$ and\\
$i \in Y \Rightarrow f_i \restriction G_{n+1} = f$.  By renaming
$\{i:i < \mu\} \subseteq Y$, now $\langle f_i f^{-1}_0:i < \mu \rangle$ is
a $\mu$-witness by part (1). 

\item  Follows by \ref{1.11}(2) and \ref{1.10}(2). 

\item  By \ref{1.11}(3) by the well foundedness of the ordinals (as in the
proof of \ref{1.10}(2),(8). 

\item  Because for $i < j,(f_jf^{-1}_i)^{-1} = (f_i f^{-1}_j)$ and
\ref{1.7}(2).  
\end{enumerate}
\end{proof}

\begin{convention}\label{1.12}
By renaming and \ref{1.11}(4), \wilog \,
  $\bar \alpha^* \in \Delta_n$ for every $n$.
\end{convention}
\bigskip

\begin{claim}\label{1.13}  Each $\alpha^*_t(t \in \mathbf I)$ is a 
non-successor ordinal (i.e. limit or zero).
\end{claim}

\begin{proof}  Fix $n < \omega$.

Assume $s \in \mathbf I$ is a counterexample.  So $\alpha^*_s = \beta^* +1,
\beta^* \ge 0$.  Let $\bar \beta = \langle \beta_t:t \in \mathbf I \rangle$
be defined as follows: $\beta_t$ is $\alpha_t$ if $t \ne s$ and is $\beta^*$ 
if $t=s$.  We shall prove that $\bar \beta \in \Gamma_{n+1}$ thus getting a
contradiction.  So let $\mu < \lambda$ and we shall find a $\mu$-witness for
$\bar \beta \in \Gamma_{n+1}$.  Let $\mu'$ be such that $\mu 
|G_{n+1}| < \mu' < \lambda$.  As $\bar \alpha^* \in \Gamma_n$ (see
\ref{1.12}) there is a $\mu'$-witness $\langle f_i:i < \mu' \rangle$ for
$\bar \alpha^* \in \Gamma_n$, as earlier \wilog \, $i < \mu \Rightarrow f_i
\restriction G_{n+1} = f$ for some $f$.  We shall prove that $\langle f_i
f^{-1}_0:i < \mu \rangle$ is a $\mu$-witness for $\bar \beta \in 
\Gamma_{n+1}$.   Let $f'_i = f_if^{-1}_0$ for $i < \mu$.

Clause (a):

$f' \restriction G_{n+1} = (f_0f^{-1}_0) \restriction G_{n+1} = 
e_{G_{n+1}}$ because $f_i \restriction G_{n+1} = f_0 \restriction G_{n+1}$.

Clause (b):

Let $i \ne j < \mu$.  If $t \in \mathbf I \backslash \{s\}$ then\\ 
rk$_t(e_{G_{n+1}},f'_i(f'_j)^{-1}) = \text{ rk}_t(e_{G_{n+1}},f_if^{-1}_j)
\le \text{ rk}_t(e_{G_n},f_if^{-1}_j) \le \alpha^*_t = \beta_t$. 
(Why?  By group theory, by \ref{1.6}(3)$(\alpha)$, by choice of $\bar f$,
by choice of $\beta_t$, respectively).

If $t=s$, then rk$_t(e_{G_n},f_if^{-1}_j) < \text{ rk}_t(e_{G_{n+1}},
f_if^{-1}_j)$ by \ref{1.6}(3)$(\beta)$, and proceed as above.
\end{proof}

\begin{notation}\label{1.14}  For $\alpha \le \omega$ let $T_\alpha := 
\prod_{k < \alpha}\lambda_k,T := \prod_{n < \omega} T_n$ (note: treeness
used).
\end{notation}

\begin{claim} \label{1.15} There are for $n < \omega$, a sequence $\langle
f_{n,i}:i < \lambda_n \rangle$ and an ordinal $\gamma^t_n <
\alpha^*_t$ ($\alpha^*_t$ is the ordinal from \ref{1.12}) such that

\begin{enumerate}

\item  $f_{n,i} \in G_\omega,f_{n,i} \restriction G_{n+1} = 
e_{G_{n+1}}$ for all $i< \lambda_n$;

\item   for each $t \in \mathbf I$ for every $h \in H^t_n$ and
$i < j < \lambda_n$ we have: \\
rk$_t(h,f_{n,i}f^{-1}_{n,j}) \le \gamma^t_n$;

\item   rk$_t(e_{H^t_n},f_{n,i}f^{-1}_{n,j}) \ge 
\gamma^t_{n-1}$ for $i < j < \lambda_n$\\ and $\gamma^t_{n-1} 
\ge 0 \Rightarrow \text{ rk}_t(e_{H^t_n},f_{n,i}f^{-1}_{n,j}) > 
\gamma^t_{n-1}$

\item   $\gamma^t_{n-1} < \gamma^t_n$ if $\alpha^*_t > 0$ and
$\gamma^t_n =-1$ if $\alpha^*_t = 0$.

\end{enumerate}
\end{claim}
\bigskip

We delay the rest of proof  for a while.

\begin{convention}  Let $\gamma^t_n,g_{n,i} \, (n < \omega,
i < \lambda_n)$ be as in \ref{1.15}.

\end{convention}
\bigskip

\begin{definition}\label{1.17}  We set $f_\eta = 
g_{n-1,\eta(n-1)}g_{n-2,\eta(n-2)} \ldots g_{0,\eta(0)}$ for $\eta \in
T_n$.  Then define $f_\eta$ for $\eta \in T_\omega$ as follows: 
$f_\eta$ is the element of $G^\omega$ satisfying $f_\eta \restriction G_n = 
f_{\eta \restriction n}$.  It is well defined by:
\end{definition}

\begin{fact}\label{1.18}  \begin{enumerate}

\item For $\eta\in T_\omega$ and $m \le n < \omega$ we
have
\[
f_{\eta \restriction n} \restriction G_{n+1} = f_{\eta \restriction m}
\restriction G_{n+1}.
\]

\item For $\eta \in T_\omega$ we have $f_\eta \in G_\omega$ is well defined (as
the inverse limit of $\langle f_{\eta \restriction n} \restriction G_n:
n < \omega \rangle$, so $n < \omega \rightarrow f_\eta \restriction G_n =
f_{\eta \restriction n}$.
\end{enumerate}
\end{fact}

\begin{proof} \begin{enumerate}

\item As $\pi_{n,\omega}$ is a homomorphism it is enough to prove\\
$(f_{\eta \restriction n}(f_{\eta \restriction m})^{-1}) \restriction
G_{n+1} = e_{G_{n+1}}$, hence it is enough to prove \\$n \le k < \omega 
\Rightarrow (f_{\eta \restriction k} f^{-1}_{\eta \restriction (k+1)}) 
\restriction G_{n+1} = e_{G_{n+1}}$ which follows from 
$k < \omega \Rightarrow f_{\eta \restriction k} 
f^{-1}_{\eta \restriction (k+1)} \restriction G_{k+1} = e_{G_{k+1}}$, 
which means\\ $f_{k,\eta(k)} \restriction G_{k+1} = e_{G_{k+1}}$ which 
holds by clause (a) of \ref{1.12}. 

\item Follows by part (1) and $G_{\omega}$ being an inverse limit.
\end{enumerate}
\end{proof}

\begin{proposition}\label{1.19}  Let $\eta,\nu \in T_\omega$.  If $\eta \ne
\nu$ and $t \in \mathbf I$, then $f_\eta f^{-1}_\nu \notin \sigma^t_\omega
(H^t_\omega)$.
\end{proposition}

\begin{proof}  Suppose for the sake of
 contradiction that for some $g \in H^t_\omega$
we have $\sigma^t_\omega(g) = f_\eta f^{-1}_\nu$.  

Let $k$ be minimal
such that
$\eta \restriction k = \nu \restriction k,\eta(k) \ne \nu(k)$, without
loss of generality 
$\eta(k) < \nu(k)$.  For $\ell \ge k$ let $\xi^\ell$ be
rk$_t(g \restriction H^t_\ell,f_{\eta \restriction (\ell + 1)}
f^{-1}_{\nu \restriction (\ell +1)})$.  We will reach a contradiction by
showing that $\ell \ge k \Rightarrow 0 \le \xi^\ell \le \gamma^t_k$ and
$\ell > k \Rightarrow \xi^{\ell +1} < \xi^\ell$.

Note

\item {$(*)_1$}  if $\ell \le \alpha \le \omega$, then
rk$_t(g \restriction H^t_\ell,f_{\eta \restriction \alpha} 
f^{-1}_{\nu \restriction \alpha}) \ge 0$ as 
$\sigma^t_\ell(g \restriction H^t_\ell) = \sigma^t(g) \restriction
G^t_\ell = (f_\eta f^{-1}_\nu) \restriction F^t_\ell$ 
and \ref{1.18}.

For $\ell =k$, we show that $\xi^k \le \gamma^t_k$.  Let $i =\eta[k],j = \nu
[k]$.  By the choice of $k,i \ne j$.  In this case 
$f_{\eta \restriction (\ell +1)}f^{-1}_{\nu \restriction (\ell +1)} =
f_{k,\eta(k)}f^{-1}_{k,\nu(k)}$ by the minimality of 
$k$ and, of course, $f_{k,\eta(k)}f^{-1}_{k,\nu(k)} = f_{k,i}f^{-1}_{k,j}$,
hence $\xi^k = \text{ rk}_t(g \restriction H^h_k,f_{k,i}f^{-1}_{k,j}) 
\le \gamma_k$ by clause (b) of \ref{1.15}. 
Note: if $\alpha^*_t = 0$, then $\gamma^t_m = -1$ for $m < \omega$ hence
$\xi^k = -1$, but $(f_\eta f^{-1}_\nu) \restriction G_k = 
(f_{\eta \restriction (k+1)}f^{-1}_{\nu \restriction (k+1)}) \restriction
G_k$ immediate contradiction.  So assume $\alpha^*_t \ge 0$ hence
$0 \le \gamma^t_m < \gamma^t_{m+1}$.

Now we proceed inductively.  We assume that $\xi^\ell \le \xi^k$ and show
that $\xi^{\ell +1} < \xi^\ell$.  Let $i = \eta[\ell +1],j = \nu[\ell +1]$, 
and let \\$\zeta = \text{ rk}_t(g \restriction H^t_{\ell +1},
f_{\eta \restriction(\ell +1)}f^{-1}_{\nu \restriction (\ell +1)})$.  Observe:

\item  $(*)_2$   $\zeta < \text{ rk}_t(g \restriction H^t_\ell,
f_{\eta \restriction (\ell +1)}f^{-1}_{\nu \restriction (\ell +1)}) = 
\xi^\ell$ 
[why?  by \ref{1.6}(3) and $(*)_1$ above.]

So

\item {$(*)_3$}  $\xi^{\ell +1} = \text{ rk}_t(g \restriction H^t_{\ell +1},
f_{\eta \restriction (\ell +2)} f^{-1}_{\nu \restriction (\ell +2)})$ 

$\quad = \text{ rk}_t
(g \restriction H^t_{\ell +1},f_{\ell +1,\eta(\ell +1)}
(f_{\eta \restriction (\ell +1)}f^{-1}_{\nu \restriction (\ell +1)})
f_{\ell +1,\nu(\ell +1)})$ 

$\quad = \text{ rk}_t(e_{H^t_{\ell+1}}(g \restriction H^t_{\ell +1})
e_{H^t_{\ell +1}},f_{\ell +1,\eta(\ell +1)}(f_{\eta \restriction (\ell +1)}
f^{-1}_{\nu \restriction (\ell +1)})f_{\ell + 1,\nu(\ell +1)})$.

Now:

\item $(*)_4$  rk$_t(e_{H^t_{\ell +1}},f_{\ell +1,\eta(\ell +1)}) >
\gamma^t_\ell$ 
(why?  by clause (c) of \ref{1.15})

\item $(*)_5$  rk$_t(g \restriction H^t_{\ell +1},f_{\eta \restriction
(\ell +1)}f^{-1}_{\nu \restriction (\ell +1)}) = \xi^\ell \le \xi^k \le
\gamma^t_k \le \gamma^t_\ell$ \\
(why?  the equality by the definition of $\xi^\ell$, the first inequality
by the induction hypothesis and the second inequality was proved above
(for $\ell =k$), the last inequality by \ref{1.15} clause (d)

\item $(*)_6$  rk$_t(e_{H^t_{\ell +1}},g_{\ell +1,\nu(\ell +1)}) >
\gamma^t_\ell$ 
(why?  by clause (c) of \ref{1.15}).

Hence by \ref{1.6}(3)

\item $(*)_7$  rk$_t(e_{H^t_{\ell +1}}(g \restriction H^t_{\ell +1})
e_{H^t_{\ell +1}},f_{\ell +1,\eta(\ell)}(f_{\eta \restriction (\ell +1)}
f^{-1}_{\nu \restriction (\ell +1)})f_{\ell +1,\nu(\ell +1)})$ \\
$ = \text{ rk}(g \restriction G_{\ell +1},f_{\eta \restriction (\ell +1)}
f^{-1}_{\nu \restriction (\ell +1)})$.

Together we get the induction demand for $\ell +1$. 
\end{proof}

Before proving \ref{1.15} and finishing we prove

\begin{claim}\label{1.20}  Assume $-1 \le \beta_t < \alpha^*_t$ for
$t \in \mathbf I$ and $n < \omega$ and $\mu < \lambda$.  Then we can
find $\langle f_i:i < \mu \rangle$ such that

\begin{enumerate}

\item   $f_i \in G_\omega$ and $f_i \restriction G_{n+1} =
e_{G_{n+1}}$

\item   $t \in \mathbf I $ and $i \ne j \Rightarrow \text{ rk}_t
(e_{H^t_n},f_i f^{-1}_j) \in [\beta^t,\alpha^*_t)$

\item   $t \in \mathbf I$ and $i < \mu = \text{ rk}_t(e_{H^t_n},f_i) 
\in [\beta_t,\alpha^*_t)$. 

\end{enumerate}
\end{claim}

\begin{proof}  For each $s \in \mathbf I$ we define $\bar \beta^s = \langle
\beta^s_t:t \in I \rangle$ by:
\[
\beta^s_t = \cases \alpha_t &\text{ if } t \ne s \\
  \beta_t &\text{ if } t = s \endcases
\]
So $\bar \beta^s \le \bar \alpha^*,\bar \beta^s \ne \bar \alpha^*$, so as
$\bar \alpha^* \in \times_{n < \omega} \Delta_m$ necessarily $\bar \beta^s
\notin \Gamma_n$, hence for some $\mu^s < \lambda$ there is no $\mu^s$-witness
for $\bar \beta^s$ and $n$ (check the definition of $\Gamma_n$).

Let $\mu_1 < \lambda$ be $> \mu + \text{ max}\{\mu^s:s \in \mathbf I\}$.

Let $\chi < \lambda$ be large enough (so that
it will be possible to use the 
finite Ramsey theorem  when $\lambda = \aleph_0$ and when $\lambda>\aleph_0$
the  Erd\"os Rado 
theorem  we require that  $\chi \rightarrow (\mu_1)^2_\theta$ where
$\theta = 2^{  \sum_t |H^t_n|}$).

Let $\langle f_i:i < \chi \rangle$ be a $\chi$-witness 
$\bar \alpha \in \Gamma_n$ and even $\bar \alpha \in \Gamma_{n+1}$.
For each $t \in \mathbf I,h \in H^t_n$ define the two place function
$F_{t,h}$ from $[\chi]^2$ to $\{0,1\}$ for $i<j<\chi$ let 
\[
F_{t,h}\{i,j\}:=\cases 0 &\text{ if   \quad rk}_t(h,f_i f^{-1}_j) < \beta_t\\
1 &\text{ Otherwise.}\endcases
\]

Define the two-place function $F$ from $[\chi]^2$:
For $i < j < \chi$  let $F\{i,j\} = \langle F_{t,h}(i,j):t \in \mathbf I,
h \in H^t_n \rangle$.

Clearly $|\text{Rang}(F)| \le  2^{  \sum_t |H^t_n|}$.
\medskip

Hence an application of one of the 
above partition theorems provides us with
 a set $Y \subseteq \chi,|Y| = \mu_1$ such that $F \restriction
[Y]^2$ is constant.  Without loss of generality $Y = \mu_1$.

For each $s \in \mathbf I$, clearly $\langle f_if^{-1}_0:i < \mu^s \rangle$ 
is not a $\mu^s$-witness for $\bar \beta^s$, but the only thing that may go 
wrong is the inequality, $i < j < \mu^s \Rightarrow \text{ rk}_s(e_{H^s_n},
f_i f^{-1}_j) < \beta_s$, so for some $i < j < \mu^s$ we have that
rk$_s(e_{H^s_n},f_i f^{-1}_j)
\ge \beta_s$ holds, hence

\item $(*)$  $s \in \mathbf I$ and $i < j < \mu_1 \Rightarrow
\text{ rk}_s(e_{H^s_n},f_i f^{-1}_j) \ge \beta_s$.

This means clause (b) holds and clause (a) by definition of $\langle f_i:i <
\chi \rangle$ is a $\chi$-witness for $\bar \alpha \in \Gamma_n$.  Clause
(c) follows.
So $\langle f_i:i < \mu \rangle$ is as required.
\end{proof}

\begin{proof} of \ref{1.15}  

Stipulate $\gamma^t_{-1}$: if $\alpha^*_t > 0$ it is $0$, otherwise is it
$-1$. Assume $n < \omega$ and 
$\langle \gamma^t_{n-1}:t \in \mathbf I \rangle$
is well defined, $\gamma^t_{n-1} < \alpha^*_t$.  Let $\gamma^{t,*}_n$ be:
$\gamma^t_{n-1} +1$ if $\alpha^*_t$ is a limit ordinal and $\gamma^t_{n-1}
=-1$ otherwise (i.e. $\alpha^*_t = 0$, see \ref{1.13}). 
Note that to construct the family $\{f_{n,i}:i < \lambda_n\}$ we will 
combine Claim \ref{1.20} with a second
application of the Erd\"os Rado Theorem.

Let $\theta = (2^{|H^t_n| \times |H^t_n|}) \times |\mathbf I|$ and $\chi <
\lambda$ be such that $\chi \rightarrow (\lambda_n+2)^3_\theta$ (exists by
Ramsey theorem if $\lambda = \aleph_0$ and by Erd\"os Rado theorem if
$\lambda > \aleph_0$).  Apply Claim \ref{1.20} to get a family
$\{f_i:i < \chi\}$ satisfying:

\begin{enumerate}
\item   $f_i \restriction G_{n+1} = e_{G_{n+1}}$,

\item   for $i \ne j$ and $t \in \mathbf I$, we have $\gamma^{t,*}
_{n-1} \le \text{ rk}_t(e_{H^t_n},f_if_j^{-1}) < \alpha^*_t$.
\end{enumerate}

For $t \in \mathbf I,\bar g = \langle g_1,g_2 \rangle,g_1,g_2 \in H^t_n$
such that $\sigma^t_n(g) = e_{G_n}$ define a coloring $F_{t,\bar g}$ of
$[I]^3$ by two colors according to the following scheme:  for
$\varepsilon < \zeta < \xi < \chi$, let 
\[ F_{t,g}\{\varepsilon,\zeta,\xi\}:=\cases
\text{\it{red}} \qquad \text{ if } \quad \text{ rk}_t(g_1,f_{i_\varepsilon}
f^{-1}_\zeta) \le \text{ rk}_t(g_2,f_\zeta f^{-1}_\xi);
\\
\text{\it{green}} \qquad \text{ if } \quad \text{ rk}_t(g_1,f_{i_\varepsilon}
f^{-1}_\zeta) > \text{ rk}_t(g_2,f_\zeta f^{-1}_\xi)\endcases.
\]
By the Ramsey theorem (if $\lambda = \aleph_0$) or Erd\"os Rado Theorem 
if $\lambda > \aleph_0$ there is a set $J \subseteq \chi,\text{ otp}
(J) = \lambda_n+2$ such that each coloring is constant on $[J]^3$.  
Let the value of $F_{t,\bar g}$ on $[J]^3$ be denoted $c_{t,\bar g}$.  
Observe that $c_{t,\bar g}$ is never {\it{green}} as this would 
produce a descending $\omega$-sequence of ordinals as if
$\varepsilon_\ell \in J,\varepsilon_\ell < \varepsilon_{\ell+1}$ for $\ell
< \omega$, then rk$_t(g,f_{\varepsilon_\ell} f^{-1}_{\varepsilon_{\ell +1}})
> \text{rk}_t(g,f_{\varepsilon_{\ell +1}} f^{-1}_{\varepsilon_{\ell+2}})$,
so $\langle \text{rk}_t(g,f_{\varepsilon_{2 \ell}} 
f^{-1}_{\varepsilon_{2 \ell +1}}):\ell < \omega \rangle$ 
is strictly decreasing.  

Let $\varepsilon(*) = \text{ Min}(J)$ and
$J_0 = \{\varepsilon \in J:\text{otp}(\varepsilon \cap J) < \lambda\}$ and
$\alpha$ is the $\lambda_n$-th member of $J,\beta$ the $(\lambda_n+1)$-th
member of $J$ and let $\gamma^t_n = \text{ rk}_t(e_{H^t_n},f_\alpha
f^{-1}_\beta)$, by clause (b) above $\gamma^{t,*}_n \le \gamma^t_n < 
\alpha^*_t$ so $\alpha^*_t = 0 \Rightarrow \gamma^t_n = -1$ and 
$\alpha^*_t > 0 \Rightarrow \gamma^t_{n-1} < \gamma^t_n$.

We claim that $\{f_if^{-1}_{\varepsilon(*)}:i \in J_0\}$ \, 
(remember $J_0 \subseteq J,|J_0| = \lambda_n)$ provides a set that can 
play the role of $\{f_{n,i}:i < \lambda_n\}$.  We note

\item $(*)_1$  rk$_t(g,f_\varepsilon f^{-1}_\zeta) \le \gamma^n_t$ 
for $\varepsilon < \zeta$ in $J_0$ 
[why?  clearly $\alpha < \beta < \varepsilon < \zeta$ are in $J$ hence by the
choice of $J$ we have rk$_t(g,f_\varepsilon f^{-1}_\zeta) \le
\text{ rk}_t(g,f_\zeta f^{-1}_\alpha) \le \text{ rk}_t(g,f_\alpha
f^{-1}_\beta) = \gamma^t_n$].

Now clauses (1), (4) of \ref{1.15} holds by clause (1) above, clause (3) of
\ref{1.15} holds by $(*)_1$ and clause (4) of \ref{1.15} holds by the 
choice of the $\gamma^*_t$.  We are left with clause (2). 
Let $h \in H^t_n$, as above clearly for $\Upsilon < \xi < \zeta < \xi$ in
$J$ we have 
rk$_t(h,f_\varepsilon f^{-1}_\zeta) \le \text{ rk}_t(h,f_\zeta
f^{-1}_\xi)$.  Hence for 
$\Upsilon \varepsilon < \zeta < \xi$ in $J_0$ we have
\[
\gamma^n_t \ge \text{ rk}_t(e_{H^t_n},f_\varepsilon f^{-1}_\zeta)
\]
\[ =\text{ rk}
_t(h^{-1},(f_\varepsilon f^{-1}_\xi)(f_\zeta f^{-1}_\xi)^{-1}) 
\]
\[
 \ge 
\text{Min}\{\text{rk}_t(h,f_\varepsilon f^{-1}_\xi),\text{ rk}_t(h^{-1},
(f_\zeta f^{-1}_\xi)^{-1}\} 
\]
\[= 
\text{ Min}\{\text{rk}_t(h,f_\varepsilon
f^{-1}_\xi),\text{ rk}_t(h,f_\zeta f^{-1}_\xi)\} 
\]
\[\ge 
\text{Min}\{\text{rk}_t(h,f_\Upsilon f^{-1}_\varepsilon), 
\text{ rk}_t(h,f_\Upsilon
f^{-1}_\varepsilon)\} 
\]
\[= \text{ rk}_t(h,f_\Upsilon f^{-1}_\varepsilon).
\]
So 
giving also clause (2) of \ref{1.15}.
\end{proof}
\hfill$\square_{\ref{1.1}}$
\bigskip

\begin{remark}
The result about the cardinality of  $Ext_p(G,\mathbf {Z})$ can
be derived from Theorem \ref{1.1} using the following
definition (which constructs an isomorphic group ot  $Ext_p(G,\mathbf {Z})$).
\end{remark}

\begin{definition}\label{extp} Given an abelian group $G$, 
let $G^*:=Hom(G,\mathbf {Z})$ and for a prime $p$
denote by  $G^p$ the group $Hom(G,\mathbf {Z}/p\mathbf {Z})$.  For
$g\in G^*$ let $g\mapsto g/p$ be the natural homomorphism from $G^*$
into $G^p$.  By $G^*/p$ denote the subgroup of $G^p$ 
which is the image of$G^*/p$ under $g\mapsto g/p$. Finally
\[
Ext_p(G,\mathbf {Z}):=G^p/(G^*/p).
\]
\end{definition}

Recall that when $\lambda$ is $\aleph_0$
or strong limit of cofinality $\aleph_0$ then $\lambda^{\aleph_0}=2^{\lambda}$.

The group $H_{\omega}$ corresponde to the subgroup $G^*/p$
and the $\sigma$'s are inclusions.

We have learned from Paul Eklof that Christian U. Jensen in his book \cite{Jen}
have a proof of Theorem 1.0 of \cite{GrSh} for the case that $\lambda=\aleph_0$.

\end{document}